\documentclass[12pt]{report}
\usepackage{amssymb,amsmath,makeidx,verbatim}
\newcommand{\R}{\mathbb{R}}

\newcommand{\C}{\mathbb{C}}
\newcommand{\Z}{\mathbb{Z}}

\newcommand{\be}{\begin{enumerate}}
\newcommand{\ee}{\end{enumerate}}
\newcommand{\bq}{\begin{eqnarray*}}
\newcommand{\eq}{\end{eqnarray*}}

\begin{document}
%\pagenumbering{roman}
\newcommand{\disp}{\displaystyle}
\thispagestyle{empty}
\begin{center}
\textbf{On harmonic analysis of spherical convolutions on semisimple Lie groups.\\}
\ \\
\textbf{Olufemi O. Oyadare}\\
\ \\
{\bf Abstract}
\end{center}
\begin{quote}

\indent This paper contains a non-trivial generalization of the Harish-Chandra transforms on a connected semisimple Lie group $G,$ with finite center, into what we term spherical convolutions. Among other results we show that its integral over the collection of bounded spherical functions at the identity element $e \in G$ is a weighted Fourier transforms of the Abel transform at $0.$ Being a function on $G,$ the restriction of this integral of its spherical Fourier transforms to the positive-definite spherical functions is then shown to be (the non-zero constant multiple of) a positive-definite distribution on $G,$ which is tempered and invariant on $G=SL(2,\R).$ These results suggest the consideration of a calculus on the Schwartz algebras of spherical functions. The Plancherel measure of the spherical convolutions is also explicitly computed.
\end{quote}
\ \\
\textbf{Subject Classification:} $43A85, \;\; 22E30, \;\; 22E46$\\
\textbf{Keywords:} Spherical Bochner theorem: Tempered invariant distributions: Harish-Chandra's Schwartz algebras\\
\ \\
{\bf $\bf{1}\;\;\;\;$ Introduction}\\

Let $G$ be a connected semisimple Lie group with finite center,
and denote the Harish-Chandra-type Schwartz spaces of functions on $G$
by ${\cal C}^p(G)$, $0 < p \leq 2.$ We know that ${\cal C}^p(G)\subset L^p(G)$ for every such $p$, and if $K$ is a maximal
compact subgroup of $G$ such that ${\cal C}^p(G//K)$ represents the subspace
of ${\cal C}^p(G)$ consisting of the $K-$bi-invariant functions, Trombi and Varadarajan [$11.$] have shown that the spherical Fourier
transform $f \mapsto \widehat{f}$ is a linear topological isomorphism
of ${\cal C}^p(G//K)$ onto the spaces $\bar{\mathcal{Z}}({\mathfrak{F}}^{\epsilon})$,
$\epsilon = \left(2/p\right)-1,$
\ \\
\ \\
$\overline{Department\;of\;Mathematics,}$
Obafemi Awolowo University,
Ile-Ife, $220005,$ Nigeria.
\text{E-mail: \textit{femi\_oya@yahoo.com}}
\ \\
consisting of rapidly decreasing functions on certain sets
${\mathfrak{F}}^{\epsilon}$ of elementary spherical functions. It
then follows that every positive-definite distribution on
$C^{\infty}_{c}(G)$ can be uniquely extended to ${\cal C}^1(G//K)$.

Using these, and improving on the results of Godement $[6.]$ on the
\textit{Bochner theorem}, Barker [$3.$] has shown that every
positive-definite distribution, $T,$ on $G$ extends uniquely to a
continuous linear functional on ${\cal C}^1(G)$ and that $$T[f] = \int_{\mathcal{P}}\widehat{f}d\mu$$ for a uniquely defined
\textit{Borel measure}, $\mu$.  Here $f \in {\cal C}^1(G//K)$ and $\mathcal{P}$
is the space of positive-definite spherical functions on $G$. This
is his \textit{spherical Bochner theorem} which has been extended to all
${\cal C}^p(G//K)$, $1 \leq p \leq 2$, with the requirement that $supp(\mu)\subset
\mathfrak{F}^{\epsilon}.$

Now if $f \in C^{\infty}_c(G)$ and $\varphi_{\lambda} \in {\cal C}^p(G//K)$
we define a function on $G,$ termed \textit{spherical convolutions} and denoted $\mathcal{H}_{x,\lambda}f,$ as $$\mathcal{H}_{x,\lambda}f=(f\ast\varphi_{\lambda})(x).$$
We show, among other properties, that the map $\lambda \mapsto \mathcal{H}_{x,\lambda}f$ is well-defined on $\mathcal{P},$
{\it Weyl group} invariant, and that the integral over $\mathcal{P}$ of its spherical Fourier transform is a non-zero constant multiple of
$T[\varphi_{\lambda}]$ for every $f \in{\cal C}^p(G//K)$ whenever $supp(\mu)\subset \mathfrak{F}^{\epsilon}.$ This gives an expansion formula for this integral when $\varphi_{\lambda} \in \mathcal{C}^2_\tau(SL(2,\R)),$ where $\tau$ is a double representation on $K=SO(2).$ When considered as the function $x \mapsto \mathcal{H}_{x,\lambda}f$ on $G,$ the behaviours of the spherical convolutions at the identity element $x=e$ and at $\lambda=0$ show both its generalization of the \textit{Harish-Chandra transforms} and its relationship with the \textit{elementary} spherical function $\Xi$ respectively. Its membership of  the Schwartz algebra ${\cal C}^{2}(G//K),$ which leads to the consideration of its spherical Fourier transforms, results to the proof of a more inclusive \textit{Plancherel formula} for ${\cal C}^{2}(G//K).$

Details of these results are contained in $\S4.$ after giving a preliminary on the structure theory of $G$ in $\S2.$ and the spherical Bochner theorems in $\S3.$\\
\ \\
{\bf $\bf{2}\;\;\;\;$ Preliminaries}\\

For the connected semisimple Lie group
$G$ with finite center, we denote its Lie algebra by $\mathfrak{g}$
whose \textit{Cartan decomposition} is given as $\mathfrak{g} = \mathfrak{t}\oplus\mathfrak{p}.$ Denote by $\theta$ the \textit{Cartan involution} on $\mathfrak{g}$ whose collection of fixed points is $\mathfrak{t}.$
We also denote by $K$ the analytic subgroup of $G$ with Lie
algebra $\mathfrak{t}.$  $K$ is then a maximal compact subgroup of $G.$
Choose a maximal abelian subspace  $\mathfrak{a}$ of $\mathfrak{p}$ with algebraic
dual $\mathfrak{a}^*$ and set $A =\exp \mathfrak{a}.$  For every $\lambda \in \mathfrak{a}^*$ put
$$\mathfrak{g}_{\lambda} = \{X \in \mathfrak{g}: [H, X] =
\lambda(H)X, \forall  H \in \mathfrak{a}\},$$ and call $\lambda$ a restricted
root of $(\mathfrak{g},\mathfrak{a})$ whenever $\mathfrak{g}_{\lambda}\neq\{0\}$.
Denote by $\mathfrak{a}'$ the open subset of $\mathfrak{a}$
where all restricted roots are $\neq 0,$ and call its connected
components the \textit{Weyl chambers.}  Let $\mathfrak{a}^+$ be one of the Weyl
chambers, define the restricted root $\lambda$ positive whenever it
is positive on $\mathfrak{a}^+$ and denote by $\triangle^+$ the set of all
restricted positive roots. Members of $\triangle^+$ which form a basis for $\triangle$ and can not be written as a linear combination of other members of $\triangle^+$ are called \textit{simple.} We then have the \textit{Iwasawa
decomposition} $G = KAN$, where $N$ is the analytic subgroup of $G$
corresponding to $\mathfrak{n} = \sum_{\lambda \in \triangle^+} \mathfrak{g}_{\lambda}$,
and the \textit{polar decomposition} $G = K\cdot
cl(A^+)\cdot K,$ with $A^+ = \exp \mathfrak{a}^+,$ and $cl(A^{+})$ denoting the closure of $A^{+}.$

If we set $M = \{k
\in K: Ad(k)H = H$, $H\in \mathfrak{a}\}$ and $M' = \{k
\in K : Ad(k)\mathfrak{a} \subset \mathfrak{a}\}$ and call them the
\textit{centralizer} and \textit{normalizer} of $\mathfrak{a}$ in $K,$ respectively, then (see $[7.]$, p. $284$);
(i) $M$ and $M'$ are compact and have the same Lie algebra and
(ii) the factor  $\mathfrak{w} = M'/M$ is a finite group called the \textit{Weyl
group}. $\mathfrak{w}$ acts on $\mathfrak{a}^*_{\C}$ as a group of linear
transformations by the requirement $$(s\lambda)(H) =
\lambda(s^{-1}H),$$ $H \in \mathfrak{a}$, $s \in \mathfrak{w}$, $\lambda \in
\mathfrak{a}^*_\mathbb{\C}$, the complexification of $\mathfrak{a}^*$.  We then have the
\textit{Bruhat decomposition} $$G = \bigsqcup_{s\in \mathfrak{w}} B m_sB$$ where
$B = MAN$ is a closed subgroup of $G$ and $m_s \in M'$ is the
representative of $s$ (i.e., $s = m_sM$). The Weyl group invariant members of a space shall be denoted by the superscript $^{\mathfrak{w}}.$

Some of the most important functions on $G$ are the \textit{spherical
functions} which we now discuss as follows.  A non-zero continuous
function $\varphi$ on $G$ shall be called a \textit{(zonal) spherical
function} whenever $\varphi(e)=1,$ $\varphi \in C(G//K):=\{g\in
C(G)$: $g(k_1 x k_2) = g(x)$, $k_1,k_2 \in K$, $x \in G\}$ and $f*\varphi
= (f*\varphi)(e)\cdot \varphi$ for every $f \in C_c(G//K),$ where $(f \ast g)(x):=\int_{G}f(y)g(y^{-1}x)dy.$  This
leads to the existence of a homomorphism $\lambda :
C_c(G//K)\rightarrow \C$ given as $\lambda(f) = (f*\varphi)(e)$.
This definition is equivalent to the satisfaction of the functional relation $$\int_K\varphi(xky)dk = \varphi(x)\varphi(y),\;\;\;x,y\in G.$$

It has been shown by Harish-Chandra [$8.$] that spherical functions on $G$
can be parametrized by members of $\mathfrak{a}^*_{\C}.$  Indeed every
spherical function on $G$ is of the form $$\varphi_{\lambda}(x) = \int_Ke^{(i\lambda-p)H(xk)}dk,\; \lambda
\in \mathfrak{a}^*_{\C},$$  $\rho =
\frac{1}{2}\sum_{\lambda\in\triangle^+} m_{\lambda}\cdot\lambda,$ where
$m_{\lambda}=dim (\mathfrak{g}_\lambda),$ and that $\varphi_{\lambda} =
\varphi_{\mu}$ iff $\lambda = s\mu$ for some $s \in \mathfrak{w}.$ Some of
the well-known properties of spherical functions are $\varphi_{-\lambda}(x^{-1}) =
\varphi_{\lambda}(x),$ $\varphi_{-\lambda}(x) =
\bar{\varphi}_{\bar{\lambda}}(x),$ $\mid \varphi_{\lambda}(x) \mid\leq \varphi_{\Re\lambda}(x),$ $\mid \varphi_{\lambda}(x)\mid\leq \varphi_{i\Im\lambda}(x),$ $\varphi_{-i\rho}(x)=1,$ $\lambda \in \mathfrak{a}^*_{\C},$ while $\mid \varphi_{\lambda}(x) \mid\leq \varphi_{0}(x),\;\lambda \in i\mathfrak{a}^{*},\;x \in G.$ Also if $\Omega$ is the \textit{Casimir operator} on $G$ then
$$\Omega\varphi_{\lambda} = -(\langle\lambda,\lambda\rangle +
\langle \rho, \rho\rangle)\varphi_{\lambda},$$ where $\lambda \in
\mathfrak{a}^*_{\C}$ and $\langle\lambda,\mu\rangle
:=tr(adH_{\lambda} \ adH_{\mu})$ for elements $H_{\lambda}$, $H_{\mu}
\in {\mathfrak{a}}.$ The elements $H_{\lambda}$, $H_{\mu}
\in {\mathfrak{a}}$  are uniquely defined by the requirement that $\lambda
(H)=tr(adH \ adH_{\lambda})$ and $\mu
(H)=tr(adH \ adH_{\mu})$ for every $H \in {\mathfrak{a}}$ ([$7.$],
Theorem $4.2$). Clearly $\Omega\varphi_0 = 0.$

Due to a hint dropped
by Dixmier $[5.]$ $(cf.\;[10.])$ in his discussion of some functional calculus,
it is necessary to recall the notion of
a \textit{`positive-definite'} function and then discuss the situation for
positive-definite spherical functions.  We call a continuous function
$f: G \rightarrow \C$ (algebraically) positive-definite whenever, for all
$ x_1,\dots,x_m $ in $G$ and all $ \alpha_1,\dots,\alpha_m$ in $\C,$ we have $$\sum^m_{i,j=1}\alpha_i\bar{\alpha}_jf(x^{-1}_i x_j) \geq 0.$$  It
can be shown $(cf.\;[7.])$ that $f(e) \geq 0$ and $|f(x)| \leq f(e)$ for every
$x \in G$ implying that the space ${\cal P}$ of all
positive-definite spherical functions on $G$ is a subset of the
space ${\mathfrak{F}}^{1}$ of all bounded spherical functions on $G.$

We know, by the Helgason-Johnson theorem ($[9.]$), that $${\mathfrak{F}}^{1}=
\mathfrak{a}^*+iC_{\rho}$$ where $C_{\rho}$ is the convex hull of $\{s\rho: s \in
\mathfrak{w}\}$ in $\mathfrak{a}^*.$ Defining the \textit{involution} $f^*$ of $f$ as $f^*(x) =
\overline{f(x^{-1})}$, it follows that $f = f^*$ for every $f \in
{\cal P}$, and if $\varphi_{\lambda} \in {\cal P}$, then $\lambda$
and $\bar{\lambda}$ are Weyl group conjugate, leading to a realization of $\mathcal{P}$ as a subset of $\mathfrak{w} \setminus \mathfrak{a}^*_{\C}.$  ${\cal P}$ becomes
a locally compact Hausdorff space when endowed with the \textit{weak $^{*}-$topology} as a subset of $L^{\infty}(G)$.\\
\ \\
{\bf $\bf{3}\;\;\;\;$ The Spherical Bochner Theorem and its Extension}\\

Let $$\varphi_0(x):= \int_{K}\exp(-\rho(H(xk)))dk$$ be denoted
as $\Xi(x)$ and define $\sigma: G \rightarrow \C$ as
$$\sigma(x) = \|X\|$$ for every $x = k\exp X \in G,\;\; k \in K,\; X
\in \mathfrak{a},$ where $\|\cdot\|$ is a norm on the finite-dimensional
space $\mathfrak{a}.$ These two functions are spherical functions on
$G$ and there exist numbers $c,d$ such that $$1 \leq \Xi(a)
e^{\rho(\log a)} \leq c(1+\sigma(a))^d.$$ Also there exists $r_0
> 0$ such that $c_0 =: \int_G\Xi(x)^2(1+\sigma(x))^{r_0}dx
< \infty$ ($[13.],$ p. $231$).  For each
$0 \leq p \leq 2$ define ${\cal C}^p(G)$ to be the set consisting of
functions $f$ in $C^{\infty}(G)$ for which $$\|f\|_{g_1,
g_2;m} :=\sup_G|f(g_1; x ; g_2)|\Xi (x)^{-2/p}(1+\sigma(x))^m <
\infty$$ where $g_1,g_2 \in \mathfrak{U}(\mathfrak{g}_{\C}),$ the \textit{universal
enveloping algebra} of $\mathfrak{g}_{\C},$ $m \in \Z^+, x \in G,$
$f(x;g_2) := \left.\frac{d}{dt}\right|_{t=0}f(x\cdot(\exp tg_2))$
and $f(g_1;x) :=\left.\frac{d}{dt}\right|_{t=0}f((\exp
tg_1)\cdot x).$ We call ${\cal C}^p(G)$ the Schwartz space on $G$
for each $0 < p \leq 2$ and note that ${\cal C}^2(G)$ is the
well-known (see $[1.]$) Harish-Chandra space of rapidly decreasing functions on
$G.$ The inclusions $$C^{\infty}_{c}(G) \subset {\cal C}^p(G)
\subset L^p(G)$$ hold and with dense images. It also follows that
${\cal C}^p(G) \subseteq {\cal C}^q(G)$ whenever $0 \leq p \leq q
\leq 2.$ Each ${\cal C}^p(G)$ is closed under \textit{involution} and the
\textit{convolution}, $*.$ Indeed ${\cal C}^p(G)$ is a Fr$\acute{e}$chet algebra ($[12.],$ p. $69$). We endow ${\cal C}^p(G//K)$
with the relative topology as a subset of ${\cal C}^p(G)$.

For any measurable function $f$ on $G$ we define the \textit{spherical Fourier
transform} $\widehat{f}$ as $$\widehat{f}(\lambda) = \int_G f(x)
\varphi_{-\lambda}(x)dx,$$ $\lambda \in \mathfrak{a}^*_{\C}.$ It
is known (see $[3.]$) that for $f,g \in L^1(G)$ we have:\\
\begin{enumerate}
\item [$(i.)$] $(f*g)^{\wedge} = \widehat{f}\cdot\widehat{g}$ on $ {\mathfrak{F}}^{1}$
whenever $f$ (or $g$) is right - (or left-) $K$-invariant; \item
[$(ii.)$] $(f^*)^{\wedge}(\varphi) =
\overline{\widehat{f}(\varphi^*)}, \varphi \in {\mathfrak{F}}^{1}$; hence
$(f^*)^{\wedge} = \overline{\widehat{f}}$ on ${\cal P}:$ and, if we
define $f^{\#}(g) := \int_{K\times K}f(k_1xk_2)dk_1dk_2,  x\in
G,$ then \item [$(iii.)$] $(f^{\#})^{\wedge}=\widehat{f}$ on ${\mathfrak{F}}^{1}.$
\end{enumerate}

In order to know the image of the spherical Fourier transform when
restricted to ${\cal C}^p(G//K)$ we need the following spaces that are central to the statement
of the well-known result of Trombi and Varadarajan [$11.$] (Theorem $3.1$ below).

Let $C_\rho$ be the closed convex hull of the (finite) set $\{s\rho :
s\in \mathfrak{w}\}$ in $\mathfrak{a}^*$, i.e., $$C_\rho =
\left\{\sum^n_{i=1}\lambda_i(s_i\rho) : \lambda_i \geq 0,\;\;\sum^n_{i=1}\lambda_i = 1,\;\;s_i \in \mathfrak{w}\right\}$$ where we recall that, for every
$H \in \mathfrak{a},$ $$(s\rho)(H) = \frac{1}{2} \sum_{\lambda\in\triangle^+}
 m_{\lambda}\cdot\lambda (s^{-1}H).$$  Now for each
$\epsilon > 0$ set ${\mathfrak{F}}^{\epsilon} = \mathfrak{a}^*+i\epsilon
C_\rho.$ Each ${\mathfrak{F}}^{\epsilon}$ is convex in $\mathfrak{a}^*_{\C}$ and
$$int({\mathfrak{F}}^{\epsilon}) =
\bigcup_{0<\epsilon'<\epsilon}{\mathfrak{F}}^{\epsilon^{'}}$$
([$11.$], Lemma $(3.2.2)$).  Let us define $\mathcal{Z}({\mathfrak{F}}^{0}) = \mathcal{S}
(\mathfrak{a}^*)$ and, for each $\epsilon>0,$ let
$\mathcal{Z}({\mathfrak{F}}^{\epsilon})$ be the space of all $\C$-valued
functions $\Phi$ such that  $(i.)$ $\Phi$ is defined and holomorphic
on $int({\mathfrak{F}}^{\epsilon}),$ and $(ii.)$ for each holomorphic
differential operator $D$ with polynomial coefficients we have $\sup_{int({\mathfrak{F}}^{\epsilon})}|D\Phi| < \infty.$ The space
$\mathcal{Z}({\mathfrak{F}}^{\epsilon})$ is converted to a Fr$\acute{e}$chet algebra by equipping it with the
topology generated by the collection, $\| \cdot \|_{\mathcal{Z}({\mathfrak{F}}^{\epsilon})},$ of seminorms given by $\|\Phi\|_{\mathcal{Z}({\mathfrak{F}}^{\epsilon})} := \sup_{int({\mathfrak{F}}^{\epsilon})}|D\Phi|.$  It is known that $D\Phi$ above extends to
a continuous function on all of ${\mathfrak{F}}^{\epsilon}$
([$11.$], pp. $278-279$).  An appropriate subalgebra of
$\mathcal{Z}({\mathfrak{F}}^{\epsilon})$ for our purpose is the closed
subalgebra $\bar{\mathcal{Z}}({\mathfrak{F}}^{\epsilon})$ consisting of
$\mathfrak{w}$-invariant elements of $\mathcal{Z}({\mathfrak{F}}^{\epsilon})$,
$\epsilon \geq 0.$ The following well-known result affords us the
opportunity of defining a distribution on ${\cal C}^p(G//K).$

{\bf 3.1 Theorem (Trombi-Varadarajan $[11.]$).}  \textit{Let $0 < p \leq 2$ and
set $\epsilon = \left(2/p\right)-1$.  Then the
spherical Fourier transform $f \mapsto \widehat{f}$ is a linear
topological algebra isomorphism of ${\cal C}^p(G//K)$ onto $\bar{\mathcal{Z}}
({\mathfrak{F}}^{\epsilon}).\;\;\Box$}

In order to use the above theorem to state the results of
Barker $[3.],$ we require the following notions.

{\bf 3.1 Definitions.}
\begin{enumerate}
\item [$(i.)$] \textit{A distribution  $T$ on $G$ (i.e., $T \in C^{\infty}_c(G)'$) is said to be (integrally) \textit{positive-definite} (written as $T\gg 0$) whenever $$T [f * f^*] \geq 0,$$ for $f \in C^{\infty}_c(G).$}
    \item [$(ii.)$] \textit{A distribution $T$ on $G$ is called $K$-\textit{bi-invariant} whenever $T^{\#}=T$ where
$$T^{\#}[f] := T[f^{L(k_1)R(k_2)}],$$ for $f \in C^{\infty}_c(G).$}
\item [$(iii.)$] \textit{A measure $\mu$ defined on ${\cal P}$ is said to be of \textit{polynomial growth} if there exists a holomorphic polynomial $Q$ on $\mathfrak{a}^*_{\C}$ such that $\int_{\cal P}(d\mu/|Q|)<\infty.$}
    \item [$(iv.)$] \textit{The \textit{support}, $supp (\mu),$ of a regular Borel measure $\mu$ is the smallest closed set $A$ such that $\mu (B)= 0$ for all Borel sets $B$ disjoint from $A.$}
    \end{enumerate}
The following is the first of the main results of $[3.].$

{\bf 3.2 Theorem (The spherical Bochner theorem).}  \textit{Suppose $T \in C^{\infty}_c(G)'$ and $T \gg 0$.
Then $T$ extends uniquely to an element in $(\mathcal{C}^1(G))'$ and there exists a unique positive regular Borel measure $\mu$ of polynomial growth on ${\cal P}$ such that $$T[f] = \int_{\cal P}\widehat{f}d\mu,\;\; f \in {\cal C}^1(G//K).$$ The correspondence between $T$ and $\mu$ is bijective when restricted to $K-$bi-invariant distributions, in which case the formula holds for all $f \in {\cal C}^1(G).\;\;\Box$}

The second of the main results of $[3.]$ is a consequence of the Trombi-Varadarajan theorem (Theorem $3.1$ above) and is stated as follows.

{\bf 3.3 Theorem (The extension theorem).}  \textit{Suppose $T$ is a positive-definite distribution with spherical Bochner measure $\mu$.  Then $T \in ({\cal C}^p(G//K))'$ iff supp $(\mu) \subset {\mathfrak{F}}^{\epsilon}$ where $1 \leq p \leq 2$ and $\epsilon = \left(2/p\right)-1$.  In such a case $$T[f] = \int_{\cal P}\widehat{f}d\mu,\;\;f \in {\cal C}^p(G//K).\;\;\Box$$}
\ \\
\ \\
{\bf $\bf{4}\;\;\;\;$ Spherical Convolutions}\\

We start by defining the central notion of this research work.

{\bf 4.1 Definition.} \textit{Let $f$ be any measurable function on $G.$ The spherical convolution of $f$ is the measurable function, $\mathcal{H}_{x,\lambda}f,$ on $G\times \mathfrak{a}^*_{\C}$ given by the map $$(x,\lambda)\mapsto \mathcal{H}_{x,\lambda}f:=(f\ast \varphi_{\lambda})(x),$$ where $x \in G, \lambda \in \mathfrak{a}^*_{\C}.$}

We shall refer to the map $\lambda\mapsto \mathcal{H}_{x,\lambda}f$ as the \textit{spherical convolution of $f$ at $x \in G.$} The importance of Definition $4.1$ is seen from the next Lemma (especially the realization of the spherical Fourier transforms in item $(ii.)$).

{\bf 4.1 Lemma.} \textit{Let $f,f_{1}$ and $f_{2}$ be measurable functions on $G,$ whose identity element is denoted as $e.$ Then}
\begin{enumerate}
\item [$(i.)$] \textit{$\mathcal{H}_{x,\lambda}(f_{1}\pm cf_{2})=\mathcal{H}_{x,\lambda}f_{1}\pm c\mathcal{H}_{x,\lambda}f_{2},\;\;x \in G,\;c \in \C,\;\lambda \in \mathfrak{a}^*_{\C};$}
\item [$(ii.)$] \textit{$\mathcal{H}_{e,\lambda}f=\widehat{f}(\lambda),\;\;\lambda \in \mathfrak{a}^*_{\C};$}
\item [$(iii.)$] \textit{$\mathcal{H}_{x,0}1=\int_{G}\Xi(y^{-1}x)dy,\;\;x \in G.$}
\item [$(iv.)$] \textit{$\overline{\mathcal{H}_{x,-\lambda}f}=\mathcal{H}_{x,\overline{\lambda}}\overline{f},\;\;x \in G,\;\lambda \in \mathfrak{a}^*_{\C};$}
\end{enumerate}

{\bf Proof.} Items $(i.)$ and $(iv.)$ are clear. We recall that, for any measurable function $f$ on $G,$ the spherical Fourier
transform, $\widehat{f},$ of $f$ is given as $\widehat{f}(\lambda) = \int_G f(x)
\varphi_{-\lambda}(x)dx,$ $\lambda \in \mathfrak{a}^*_{\C}.$ Since $\varphi_{-\lambda}(x)=\varphi_{\lambda}(x^{-1}),$ for every $\lambda \in \mathfrak{a}^{*}_{\C},\;x \in G,$ this may be written as $$\widehat{f}(\lambda) = \int_G f(x) \varphi_{\lambda}(x^{-1})dx=(f \ast \varphi_{\lambda})(e)=\mathcal{H}_{e,\lambda}f.$$  This proves $(ii.).$ Item $(iii.)$ follows if we recall that $\varphi_{0}(x)=\Xi(x).\;\;\Box$

Item $(iv.)$ of Lemma $4.1$ gives the \textit{functional equation} for spherical convolutions. This Lemma (especially in item $(ii.)$) explains that the harmonic analysis of $G$ has so far been explored only with the spherical convolution at $e.$ The implication of considering only $$(e,\lambda)\mapsto \mathcal{H}_{e,\lambda}f=:\widehat{f}(\lambda)$$ is that the direct contribution of the non-identity members of $G$ to its harmonic analysis are suppressed and may never be suspected or known in the context of $\widehat{f}(\lambda).$ Indeed a great deal of properties of the spherical convolutions and their contributions to harmonic analysis on $G$ would not be available if, instead of considering the entirety of the map $(x,\lambda)\mapsto \mathcal{H}_{x,\lambda}f,$ we restrict ourselves to either $\lambda \mapsto\mathcal{H}_{e,\lambda}f=\widehat{f}(\lambda)$ or $x\mapsto\mathcal{H}_{x,0}1$ or any other special case of the spherical convolutions as has been done till now.

We shall therefore show the importance of including spherical convolutions in the harmonic analysis of $G$ by giving its bounds, $\mathfrak{w}-$group transformation and differential equation. These are contained in the following Theorem while a Plancherel formula for the functions $x \mapsto \mathcal{H}_{x,\lambda}f$ on $G$ is proved after a study of its spherical Fourier transforms.

{\bf 4.1 Theorem.} \textit{Consider a measurable function $f$ on $G,\;x \in G$ and let $\lambda \in \mathfrak{a}^*_{\C}.$ Then}
\begin{enumerate}
\item [$(i.)$] \textit{$\mid \mathcal{H}_{x,\lambda}f\mid\leq \parallel f\parallel_{1},\;\;\;$ with $f \in L^{1}(G);$}
\item [$(ii.)$] \textit{$\mid \mathcal{H}_{x,\lambda}f\mid\leq \mathcal{H}_{x,\Re\lambda}f,$ $\mid \mathcal{H}_{x,\lambda}f\mid\leq \mathcal{H}_{x,i\Im\lambda}f$ and $\mid \mathcal{H}_{x,\lambda}f\mid\leq \mathcal{H}_{x,0}f,$ for $f\geq0;$}
\item [$(iii.)$] \textit{$\mathcal{H}_{x,s\lambda}f=\mathcal{H}_{x,\lambda}f,\;\;\;$ for every $s \in \mathfrak{w};$}
\item [$(iv.)$] \textit{$\Omega\mathcal{H}_{x,\lambda}f=-(\langle\lambda,\lambda\rangle +
\langle \rho, \rho\rangle)\cdot\mathcal{H}_{x,\lambda}f,\;\;\;$ for $f \in {\cal C}^{p}(G),\;0 < p \leq 2;$}
\item [$(v.)$] \textit{$\mathcal{H}_{x,-i\rho}f=\int_{G} f(y)dy.$}
\end{enumerate}

{\bf Proof.} We employ the properties of spherical functions given in $\S2$ to establish $(i.),(ii.),(iii.)$ and $(v).$ The proof of $(iv.)$ follows if we recall that $$\Omega(f\ast \varphi_{\lambda})=f\ast \Omega\varphi_{\lambda}.\;\;\Box$$

The equation established in Theorem $4.1\;(iv.),$ or any other such equation for $q \in \mathfrak{U}(\mathfrak{g}_{\C}),$ shows that the spherical convolutions inherit the differential equations satisfied by $\varphi_{\lambda}.$ This resemblance suggests the choice of the name adopted in Definition $4.1.$ It therefore has the following series expansion.

{\bf 4.1 Corollary.} \textit{The spherical convolutions, $\mathcal{H}_{h,\lambda}f$ admit the series expansion $$\mathcal{H}_{h,\lambda}f=\sum_{s\in \mathfrak{w}}c(s\lambda)\left(e^{(s\lambda-\rho)(\log h)}+\sum_{\mu\in L^+}a_{\mu}(s\lambda)e^{(s\lambda-\rho-\mu)(\log h)}\right),$$ regardless of
the functions $f \in {\cal C}^{p}(G),\;0 < p \leq 2,$ where $h \in A^{+},$ $\lambda \in\;^*\mathfrak{F}' :=\{\nu \in\;^*\mathfrak{F}:\nu\; \mbox{is regular}\},$ $L^{+}=L\setminus\{0\},$ with $$L:= \left\{\sum_{1\leq i \leq r}m_i\alpha_i : m_1,\dots, m_r\;\;\mbox{are integers}\;\geq 0\right\},$$ for the simple roots $\alpha_i, 1\leq i \leq r,$ some subset $^*\mathfrak{F}$ of $\mathfrak{a}^*_{\C}$ and coefficient functions $a_{\mu}(\lambda)$ which may be generated from the recursive relation $$\left(\langle\mu, \mu\rangle - 2\langle\mu,\lambda\rangle\right)a_{\mu}(\lambda)$$ $$=-2\sum_{\begin{array}{c}
\alpha>0,k\geq 1\\
\mu-2k\alpha\in L
\end{array}}n(\alpha)\langle \lambda-\mu+2k\alpha-\rho,\alpha\rangle a_{\mu-2k\alpha}(\lambda),\;n(\alpha):= dim(\mathfrak{g_{\alpha}}).\;\;\Box$$}

We shall now consider the map $\lambda \mapsto \mathcal{H}_{e,\lambda}f$ for its differentiability and/or integrability with respect to $\lambda$ in some specified subset, $Y,$ of $\mathfrak{a}^{*}_{\C}.$ Indeed, $\mathcal{H}_{e,\lambda}f \in C_{c}(Y),$ for every $f \in C_{c}(G)$ and any subset, $Y,$ of $\mathfrak{a}^{*}_{\C}.$ This makes its integral, $\int_{Y}\mathcal{H}_{e,\lambda}f d\mu(\lambda),$ with respect to some normalised measure, $\mu,$ on $Y,$ worthy of an indepth study. To this end we define the map $f \mapsto f\{\varphi_{\lambda}\}$ on $G$ at the identity element, $e,$ as $$f\{\varphi_{\lambda}\}(e)=\int_{Y}\mathcal{H}_{e,\lambda}f d\mu(\lambda),\;\;f \in C_{c}(G),\;\varphi_{\lambda} \in {\cal C}^p(G).$$ Before considering the generality of $f\{\varphi_{\lambda}\}(x),$ for every $x \in G,$ we state our first major result on $f\{\varphi_{\lambda}\}(e)$ which gives an important application of its integral for $Y=\mathfrak{F}^{1}.$

Define the map $a \mapsto \beta_{\mathfrak{F}^{1}}(a)$ as $\beta_{\mathfrak{F}^{1}}(a)=\int_{\mathfrak{F}^{1}}e^{\nu(\log a)}d\mu(\nu),\;a \in A$ and the \textit{$\beta_{\mathfrak{F}^{1}}-$weighted Fourier transforms,} $\widetilde{f},$ of $f \in C_{c}(A)$ at $\lambda \in \mathfrak{F}^{1}$ as $$\widetilde{f}(\lambda)=\int_{A}f(a)e^{\lambda(\log a)}d\eta(a),$$ where $d\eta(a)=\beta_{\mathfrak{F}^{1}}(a)da.$ Observe that the above weighted Fourier transforms $\widetilde{f}(\lambda)$ reduces to the classical Fourier transforms $\widehat{f}(\lambda)=\int_{A}f(a)e^{\lambda(\log a)}da$ when $\beta_{\mathfrak{F}^{1}}(a)=1,\;\forall\;a \in A.$ We shall however use this (weighted) transforms only at the identity element $0 \in \mathfrak{F}^{1}$ of the vector space $\mathfrak{F}^{1}$ $(i.e.,\;\widetilde{f}(0)=\int_{A}f(a)d\eta(a)=\int_{A}f(a)\beta_{\mathfrak{F}^{1}}(a)da),$ in the next Theorem and this, as could be seen below, may not be un-connected with the fact that $$f\{\varphi_{\lambda}\}(e)=\int_{\mathfrak{F}^{1}}(\mathcal{H}_{x,\lambda}f)_{\mid_{x=e}}d\mu(\lambda)$$ is itself an evaluation at the identity element of $G.$

{\bf 4.2 Theorem.} \textit{Let $dx=e^{2 \rho(\log a)}dk\;da\;dn,$ where $dk,\;da$ and $dn$ are Haar measures on $K,\;A,$ and $N,$
respectively, with $dk$ normalised. For every $f \in C_{c}(G//K),$ let $\mathcal{A}(f)$ denote the Abel transform of $f$ defined on $A$ as $\mathcal{A}(f)(a)=e^{\rho(\log a)}\int_{N}f(an)dn.$ Then $$f\{\varphi_{\lambda}\}(e)=\widetilde{(\mathcal{A}f)}(0),\;\;\;\lambda \in \mathfrak{F}^{1}.$$}

{\bf Proof.} We only need to prove that $$f\{\varphi_{\lambda}\}(e)=\int_{A}\mathcal{A}f(a)\beta_{\mathfrak{F}^{1}}(a)da,\;\;\lambda \in \mathfrak{F}^{1}.$$

Indeed, using the Harish-Chandra parametrisation of $\varphi_{\lambda},$ we have
\begin{eqnarray*}
f\{\varphi_{\lambda}\}(e) &=& \int_{\mathfrak{F}^{1}}(f \ast \varphi_{\lambda})(e)d\mu(\lambda)\\
&=& \int_{\mathfrak{F}^{1}}\widehat{f}(\lambda) d\mu(\lambda)\;=\;\int_{\mathfrak{F}^{1}} \int_G f(x) \varphi_{-\lambda}(x)dx d\mu(\lambda)\\
&=&\int_{\mathfrak{F}^{1}} \int_G \int_{K} f(xk) e^{(\lambda-\rho)(H(xk))} dx dk d\mu(\lambda)\\
&=& \int_{\mathfrak{F}^{1}} \int_G f(y) e^{(\lambda-\rho)(H(y))} dy d\mu(\lambda)\\
&=& \int_{\mathfrak{F}^{1}} \int_{AN} f(an) e^{(\lambda+\rho)(\log a)} da dn d\mu(\lambda)\\
&=& \int_{\mathfrak{F}^{1}} \int_{A} (\mathcal{A}f)(a) e^{\lambda(\log a)} da d\mu(\lambda),
\end{eqnarray*}
which implies our result, using Fubini's theorem$.\;\;\Box$

The last result shows the importance of $f\{\varphi_{\lambda}\}(e)$ in the harmonic analysis of $G$ and prepares the ground for the consideration of results of \textit{Paley-Wiener type.} It will soon be clear that it is sufficient to take the measure $\mu$ as the Borel measure on $\mathcal{P}$ in defining the map $a \mapsto \beta_{\mathcal{P}}(a).$ Our motivation in this direction is to consider the general map $$\lambda \mapsto \mathcal{H}_{x,\lambda}f,$$ not only for the identity element $x=e \in G,$ but for other values of $G$ as well. This leads to the definition of $f\{\varphi_{\lambda}\}(x),\;x \in G,$ as $$f\{\varphi_{\lambda}\}(x)=\int_{Y}\mathcal{H}_{x,\lambda}f d\mu(\lambda),$$ $f \in C_{c}(G),\;\varphi_{\lambda} \in {\cal C}^p(G),$ whenever the integral is absolutely convergent. We already have two candidates for the position of $Y,$ namely $\mathfrak{F}^{1}$ and $\mathcal{P}.$ Among other results, it would be important to evaluate the above measure, $\mu,$ on these candidates. In the mean time we study some of the properties of $x \mapsto f\{\varphi_{\lambda}\}(x).$

{\bf 4.3 Theorem.} \textit{Let $f \in C_{c}(G//K),\;\varphi_{\lambda} \in {\cal C}^p(G//K), 1 \leq p \leq 2,$ and $Y=\mathcal{P}.$ Then, as a function on $\mathcal{P},\;\lambda \mapsto \mathcal{H}_{x,\lambda}f$ is continuous with compact support. Indeed, $f\{\varphi_{\lambda}\} \in {\cal C}^p(G//K).$ Moreover, we have that $$f\{\varphi_{s \lambda}\}=f\{\varphi_{\lambda}\}=f\{\varphi_{s \bar{\lambda}}\},$$ for every $s \in \mathfrak{w}.$}

{\bf Proof.} The first assertion holds, since $f \in C_{c}(G//K)$ and $C_{c}(G//K)$ is dense in ${\cal C}^p(G//K).$ The properties of $\varphi_{\lambda}$ at the end of \S 2 imply the second assertion $.\;\;\Box$

We also have that $f\{\varphi_{\lambda_{1}}+c\varphi_{\lambda_{2}}\}=f\{\varphi_{\lambda_{1}}\}+cf\{\varphi_{\lambda_{2}}\},\;c \in \C,$ suggesting that the map $f \mapsto f\{\varphi_{\lambda}\}$ may be a \textit{calculus} on ${\cal C}^p(G//K).$

Let us now consider the spherical convolution map $$(x,\lambda)\mapsto \mathcal{H}_{x,\lambda}f:=(f\ast \varphi_{\lambda})(x)$$ as the function $$x \mapsto \mathcal{H}_{x,\lambda}f$$ on $G.$ Since it is measurable its spherical Fourier transforms may be computed as shown in the following result which gives how to \textit{generate} positive-definite distributions on $G$ and which will be found useful in the proof of its Plancherel formula given later in Theorem $4.7.$

{\bf 4.4 Theorem.} \textit{Let $f \in C_{c}(G//K),\;\varphi_{\lambda} \in {\cal C}^1(G//K).$ Let $\mu$ be a spherical Bochner measure corresponding to a positive-define distribution $T$ on $G.$} Then
\begin{enumerate}
\item [(i.)] \textit{$\widehat{(\mathcal{H}_{x,\lambda}f)}(\nu)=\widehat{f}(\lambda)\cdot \widehat{\varphi_{\lambda}}(\nu),\;\;x \in G,\;\nu \in \mathfrak{a}^*_{\C}.$}
\item [(ii.)] \textit{$\int_{\mathcal{P}} \widehat{(\mathcal{H}_{x,\lambda}f)}(\nu) d\mu(\nu)=\widehat{f}(\lambda)\cdot T[\varphi_{\lambda}].$}
\end{enumerate}
Moreover, if $T^{\#} = T,$ the integral in $(ii.)$ holds for all $\varphi_{\lambda} \in {\cal C}^1(G).$

{\bf Proof.} $(i.)$ Employing the defining properties of a spherical function given in $\S2.$ we have, for every $x \in G,\;\nu \in \mathfrak{a}^*_{\C},$ that $$\widehat{(\mathcal{H}_{x,\lambda}f)}(\nu)=(f\ast \varphi_{\lambda})(e)\cdot \widehat{\varphi_{\lambda}}(\nu)=\widehat{f}(\lambda)\cdot \widehat{\varphi_{\lambda}}(\nu).$$
\begin{eqnarray*}
(ii.)\;\mbox{Now fix $\varphi_{\lambda} \in C^{\infty}_c(G//K),$ then}\\ \int_{\mathcal{P}} \widehat{(\mathcal{H}_{x,\lambda}f)}(\nu) d\mu(\nu)&=& \int_{\cal P}\widehat{f}(\lambda)\cdot \widehat{\varphi_{\lambda}}(\nu)d\mu(\nu)\\
&=& \widehat{f}(\lambda)\cdot \int_{\cal P}\widehat{\varphi_{\lambda}}(\nu)d\mu(\nu)\\
&=& \widehat{f}(\lambda)\cdot T[\varphi_{\lambda}].
\end{eqnarray*}
We apply the denseness of $C^{\infty}_c(G//K)$ in ${\cal C}^1(G//K)$ to conclude the second assertion. That $(ii.)$ holds for all $\varphi_{\lambda} \in {\cal C}^1(G)$ follows from the second part of Theorem $3.2.\;\;\Box$

The extension Theorem $3.3$ leads also to an extension of Theorem $4.4$ given next.

{\bf 4.5 Theorem (Extension Theorem).}  \textit{If supp$(\mu) \subset {\mathfrak{F}}^{\epsilon}$, $\epsilon = \left(2/p\right)-1$ and $1 \leq p \leq 2$, then $\int_{\mathcal{P}} \widehat{(\mathcal{H}_{x,\lambda}f)}(\nu) d\mu(\nu)$ is a constant multiple of $T[\varphi_{\lambda}]$ for every $\varphi_{\lambda} \in {\cal C}^p(G//K).\;\;\Box$}

The conclusion on $\int_{\mathcal{P}} \widehat{(\mathcal{H}_{x,\lambda}f)}(\nu) d\mu(\nu)$ in the last Theorem above may be generalized to the Schwartz algebra ${\cal C}^p_\tau(G)$ of all $\tau-$ spherical functions on $G$ where $\tau = (\tau_1,\tau_2)$ is a \textit{double representation} of $K.$ This would be so immediately the Trombi-Varadarajan theorem, Theorem $3.1,$ is established for ${\cal C}^p_\tau(G).$ The case $p=2$ has been proved and is contained in $[1.]$ for real-rank $1$ Lie groups $G,$ and in $[2.]$ for any semisimple Lie group of any rank, while the case of general $p$ remains an open problem.

However the situation for general $p$ and the group $G= SL(2,\R),$ or its conjugate $SU(1,1),$ is contained in $[13.]$ from which other groups could be considered. Thus using the results of $[4.]$ on $\mathcal{C}^2_\tau(SL(2,\R))$ we extend the assertions of Theorem $4.5$ to all the members of $\mathcal{C}^2_\tau(SL(2,\R)).$ This leads to an expansion of $\int_{\mathcal{P}} \widehat{(\mathcal{H}_{x,\lambda}f)}(\nu) d\mu(\nu)$ for $\varphi_{\lambda}$ in the Schwartz algebras of all $\tau-$ spherical functions on $G= SL(2, \R).$ This expansion brings in the involvement of the well-known global characters of the (unitary) principal and discrete series of representations of $G=SL(2,\R)$ ($[13.]$).

To establish this expansion formula we put the needed type of measures in place. A pair $(\mu_c,\mu_d)$ is called a \textit{tempered Bochner measure pair} whenever:

$(i.)$ $\mu_c$ is a non-negative Baire measure on ${\R}$ which is symmetric and of polynomial growth. That is, $d \mu_{c}(-\lambda) =d \mu_{c}(\lambda),$ for all $\lambda \in \R$ and $$\int_{\R} d\mu_{c}(\lambda)/(1+\mid \lambda \mid^{r}) < \infty$$ for some $r \geq 0.$

$(ii.)$ $\mu_d$ is a non-negative counting measure on $ \Z^{\prime}=Z\backslash {0}$ which is of polynomial growth. That is, $$\sum_{l \in \Z^{\prime}} \mu_{d}(l)/(1+\mid l \mid^{r}) < \infty$$ for some $r \geq 0.$

The following Theorem opens up the integral contained in Theorem $4.5$ in the special case of $G=SL(2,\R)$

{\bf 4.6 Theorem (Expansion for $\int_{\mathcal{P}} \widehat{(\mathcal{H}_{x,\lambda}f)}(\nu) d\mu(\nu)$ on $\mathcal{C}^2_\tau(SL(2,\R))$).} \textit{Let $f \in C ^\infty_c(G),\;\varphi_{\lambda}\in \mathcal{C}^2_\tau(SL(2,\R)).$ Then, up to a non-zero constant, the positive-definite distribution $\int_{\mathcal{P}} \widehat{(\mathcal{H}_{x,\lambda}f)}(\nu) d\mu(\nu)$ is given as $$\int_{\mathcal{P}} \widehat{(\mathcal{H}_{x,\lambda}f)}(\nu) d\mu(\nu)= \widehat{f}(\lambda)\cdot lim_{n\rightarrow\infty}(\int^n_{-n}\Phi^\lambda[\varphi_{\lambda}]d\mu_c(\lambda)+\Sigma_{1\leq|l|\leq n}\Theta^l[\varphi_{l}]\mu_d(l)),$$ where $\Phi^\lambda$ and $\Theta^l$
are the global characters of the (unitary) principal and discrete series of representations of $G=SL(2,\R)$ and $(\mu_c,\mu_d)$ is the tempered \textit{Bochner measure} pair associated to a tempered invariant positive-definite distribution on $G.$  In particular, $\int_{\mathcal{P}} \widehat{(\mathcal{H}_{x,\lambda}f)}(\nu) d\mu(\nu)$ is a tempered invariant distribution on $G.$}

{\bf Proof.} For any tempered invariant positive-definite distribution $T$ on $G$ there corresponds a Bochner measure pair $(\mu_c,\mu_d)$
such that $$T[f]=lim_{n\rightarrow\infty}(\int^n_{-n}\Phi^\lambda[f]d\mu_c(\lambda)+\Sigma_{1\leq|l|\leq n}\Theta^l[f]\mu_d(l)).$$ This is the main result of $[4.]$ (listed there as Theorem $9.3$), which when combined with our Theorem $4.5$ gives the assertion$.\;\;\Box$

{\bf 4.1 Remark.} The expansion given above for $\int_{\mathcal{P}} \widehat{(\mathcal{H}_{x,\lambda}f)}(\nu) d\mu(\nu)$ reveals the rich structure encoded in it. Indeed since the global characters above, in terms of which it is expressed (in Theorem $4.6$), have well-known transformation under the action of the center, $\mathfrak{Z},$ of the universal enveloping algebra, $\mathfrak{U(g_{\C})},$ of the complexification $\mathfrak{g_{\C}}$ of the Lie algebra $\mathfrak{g}$ of $G,$ a study of the functional and differential equations of $\int_{\mathcal{P}} \widehat{(\mathcal{H}_{x,\lambda}f)}(\nu) d\mu(\nu)$ is very possible and suggests a harmonic analysis involving both the discrete and (unitary) principal series of, at least, $G=SL(2,\R).$

For $f \in C_{c}(G//K)$ and $\varphi_{\lambda} \in {\cal C}^p(G//K)$ as in Theorem $4.5,$ we may view the map $\lambda \mapsto f\{\varphi_{\lambda}\}$ as the \textit{evaluation} of members of $C_{c}(G//K)$ on members of ${\cal C}^p(G//K).$ This means that $\lambda \mapsto f\{\varphi_{\lambda}\}$ is an \textit{operational calculus} on the Schwartz algebras, ${\cal C}^p(G//K),$ whose spherical Fourier transform is a distribution on $G.$ This suggests the use of the term \textit{`distributional calculus'} for $f \mapsto f\{\varphi_{\lambda}\}.$ A more detailed study of $f\{\varphi_{\lambda}\}$ may therefore be conducted by considering the invariant eigendistributions on $G,$  most especially the global characters of the irreducible admissible representations of $G.$

We now consider the explicit form of the Plancherel formula for the measurable functions $x\mapsto \mathcal{H}_{x,\lambda}f$ on $G.$ A Haar measure $dx$ on $G$ is said to be \textit{admissible} if $dx=e^{2\rho (\log a)}dkdadn\;\;(x=kan)$ where $\int_{K}dk=1$ and $\int_{\overline{N}}e^{-2\rho H(\overline{n})}d\overline{n}=1$ where $d\overline{n}$ is a Haar measure on $\overline{N}:=\theta(N).$ Recall the Borel measure $d\mu(\lambda)$ from Theorem $4.4.$ The pair $(dx,d\mu(\lambda))$ of Haar measures on the pair $(G,\mathfrak{F}^{1})$ shall be de termed \textit{admissible} if, every $f$ in the Schwartz space, $\mathcal{S}(A),$ of $A$ whose Fourier transform $\widehat{f},$ already known as $$\widehat{f}(\lambda)=\int_{A}f(a)e^{\lambda(\log a)}da,\;\;\lambda \in \mathfrak{F}^{1},$$ satisfies $$f(a)=\int_{\mathfrak{F}^{1}}\widehat{f}(\lambda)e^{-\lambda(\log a)}d\mu(\lambda),\;\;a \in A.$$

{\bf 4.7 Theorem (Plancherel formula for spherical convolutions).} \textit{Let $(dy,d\mu(\lambda))$ be an admissible pair of Haar measures on the pair $(G,\mathfrak{F}^{1}),$ $x \in G$ and $f \in\mathcal{C}(G//K).$ If we define the measure $d\zeta_{x,\lambda}$ as a normalization of the spherical Bochner measure $d\mu(\lambda)$ on $\mathfrak{F}^{1}$ by the requirement that $$d\zeta_{x,\lambda}(\nu)=\frac{1}{\mid \widehat{\varphi_{\lambda}}(\nu) \mid^{2}}d\mu(\lambda),$$ then $$\int_{G}\mid f(y) \mid^{2}dy=\int_{\mathfrak{F}^{1}}\mid \widehat{(\mathcal{H}_{x,\lambda}f)}(\nu) \mid^{2}d\zeta_{x,\lambda}(\nu).$$ In particular the map $f\mapsto \widehat{(\mathcal{H}_{x,\lambda}f)},$ for $x \in G\;\mbox{and}\;\lambda \in \mathfrak{F}^{1},$ extends uniquely to a unitary isomorphism of $L^{2}(G//K)$ with $L^{2}(\mathfrak{F}^{1},d\zeta_{x,\lambda}(\nu))^{\mathfrak{w}}.$}

{\bf Proof.} Since the spherical convolutions of $f \in\mathcal{C}(G//K)$ may be considered as the functions $x \mapsto \mathcal{H}_{x,\lambda}f$ on $G$ it follows that its spherical Fourier transforms, $\nu\mapsto \widehat{(\mathcal{H}_{x,\lambda}f)}(\nu),$ is well-defined on $\mathfrak{F}^{1}.$ Therefore
\begin{eqnarray*}
\int_{\mathfrak{F}^{1}}\mid \widehat{(\mathcal{H}_{x,\lambda}f)}(\nu) \mid^{2}d\zeta_{x,\lambda}(\nu)&=& \int_{\mathfrak{F}^{1}} \mid\widehat{f}(\lambda)\mid^{2}\cdot \mid\widehat{\varphi_{\lambda}}(\nu)\mid^{2}d\zeta_{x,\lambda}(\nu)\\ \mbox{(from Theorem $4.4\;(i.)$)}\\
&=& \int_{\mathfrak{F}^{1}} \mid\widehat{f}(\lambda)\mid^{2} d\mu(\lambda)\\ \mbox{(from definition of $d\zeta_{x,\lambda}$)}\\
&=& \int_{G}\mid f(y) \mid^{2}dy\\ \mbox{(by the Plancherel formula for $f$)}.\;\;\Box
\end{eqnarray*}

The situation of Theorem $4.7$ for $x=e$ is well-known, while the inverse, $(\mathcal{H}_{x,\lambda}f)^{-1},$ for fixed $f \in\mathcal{C}(G//K),\;x \in G\;\mbox{and}\;\lambda \in \mathfrak{F}^{1},$ is given as $$(\mathcal{H}_{x,\lambda}f)^{-1}(b)(y)=\int_{\mathfrak{F}^{1}}b(\lambda)\varphi_{\lambda}(y)d\zeta_{x,\lambda}(\nu),\;b \in \mathcal{S}(\mathfrak{F}^{1})^{\mathfrak{w}},y \in G$$ and is commonly called the \textit{exact (normalized) wave packet.} The explicit expression for the Plancherel measure, $d\zeta_{x,\lambda}(\nu),$ of the spherical convolutions in terms of elementary functions of harmonic analysis is therefore given as $$d\zeta_{x,\lambda}(\nu)=\mid \mathfrak{w} \mid^{-1}\mid \widehat{\varphi_{\lambda}}(\nu) \mid^{-2}\mid c(\lambda) \mid^{-2}d\mu(\lambda),$$ for $x \in G,\;\nu,\lambda \in \mathfrak{F}^{1},$ where the map $\lambda\mapsto c(\lambda)$ is the Harish-Chandra $c-$function. The combination of Theorems $4.2$ and $4.7$ may be used to give the Paley-Wiener theorem for spherical Fourier transforms of spherical convolutions.

It is known ($[11.],\;p.\;298$) that $\mathcal{H}_{e,\lambda}f=\widehat{f}(\lambda)$ and that, in this case, the Plancherel measure, $d\zeta_{e,\lambda}(\nu),$ of the spherical convolution, $\mathcal{H}_{x,\lambda}f,$ at $x=e$ is $$d\zeta_{e,\lambda}(\nu)=\mid \mathfrak{w} \mid^{-1}\mid c(\lambda) \mid^{-2}d\mu(\lambda).$$  A non-trivial problem is to find the relation between $\nu$ and $\lambda$ for the Plancherel measure, $d\zeta_{x,\lambda}(\nu),$ of the spherical convolutions to reduce to the classical Plancherel measure, $\mid \mathfrak{w} \mid^{-1}\mid c(\lambda) \mid^{-2}d\mu(\lambda),$ on $G.$ This is equivalent to seeking those $\nu$ in terms of $\lambda$ for which $\mid \widehat{\varphi_{\lambda}}(\nu) \mid=1,$ where $\varphi_{\lambda} \in {\cal C}^1(G//K).$ We plan to address this problem in another paper.

The richness of our results, which may be ultimately seen in Theorem $4.7,$ derives from the fact that the spherical convolutions are functions on both $G$ and $\mathfrak{F}^{1}.$ This fact allows us to switch its domains between $G$ and $\mathfrak{F}^{1},$ depending on its immediate use. In all these diverse instances of the harmonic analysis on $G$ we still use the same defining functions for the spherical convolutions. We have however taken advantage of some known results in the harmonic analysis on $G$ (like the Harish-Chandra series expansion inherited by $\mathcal{H}_{h,\lambda}f$ (in Corollary $4.1$) and the classical Plancherel formula on $G$ used in the proof of Theorem $4.7$) in order to establish our results. Nevertheless our results could still be established from the scratch without recourse to the special case of $\mathcal{H}_{e,\lambda}f=\widehat{f}(\lambda).$\\

{\bf   References.}
\begin{description}
\item [{[1.]}] Arthur, J.G., \textit{Harmonic analysis of tempered distributions on semisimple Lie groups of real rank one,} Ph.D. Dissertation, Yale University, $1970.$
    \item [{[2.]}] Arthur, J.G., \textit{Harmonic analysis of the Schwartz space of a reductive Lie group,} I. II. (preprint, $1973$).
        \item [{[3.]}] Barker, W.H.,  The spherical Bochner theorem on semisimple Lie groups, \textit{J. Funct. Anal.,} vol. \textbf{20}  ($1975$), pp. $179-207.$
            \item [{[4.]}] Barker, W.H.,  Tempered invariant, positive-definite distributions on\\ $SU(1,1)/\{\pm1\},$
            \textit{Illinois J. Maths}, vol. \textbf{28}, no. $1,$ ($1984$), pp. $83-102.$
                \item [{[5.]}] Dixmier, J., Op$\acute{\mbox{e}}$rateurs de rang fini dans les repr$\acute{\mbox{e}}$sentations unitaires,\textit{ Publ. math. de l' Inst. Hautes $\acute{\mbox{E}}$tudes Scient.,} tome $\textbf{6}$ ($1960$), pp. $13-25.$
                    \item [{[6.]}] Godement, R.A.,  A theory of spherical functions, $I.$ \textit{Trans. Amer. Math. Soc.,} vol. \textbf{73} ($1952$), pp. $496-556.$
                    \item [{[7.]}] Helgason, S.,  \textit{``Differential Geometry, Lie Groups and Symmetric Spaces,''} Academic Press, New York, $1978.$
                \item [{[8.]}] Helgason, S.,  \textit{``Groups and Geometric Analysis; Integral Geometry, Invariant Differential Operators, and Spherical Functions,''} Academic Press, New York, $1984.$
                \item [{[9.]}] Helgason, S. and Johnson, K.,  The bounded spherical functions on symmetric spaces,  \textit{Advances in Math,} \textbf{3} ($1969$), pp. 586-593.
                        \item [{[10.]}] Kahane J.-P.,  Sur un th$\acute{e}$or$\grave{e}$me de Wiener-L$\acute{e}$vy, \textit{C.R. Acad. Sc., Paris,} t. $\textbf{246},$ ($1958$), pp. $1949-1951.$
                    \item [{[11.]}] Trombi, P.C. and Varadarajan, V.S.,  Spherical transforms on semisimple Lie groups, \textit{Ann. of Math.,} \textbf{94} ($1971$), pp. $246-303.$
                \item [{[12.]}] Varadarajan, V.S.,  The theory of characters and the discrete series for semisimple Lie groups, in \textit{Harmonic Analysis on Homogeneous Spaces,} (C.C.  Moore (ed.)) \textit{Proc. of Symposia in Pure Maths.,} vol. \textbf{26} ($1973$), pp. $45-99.$
                    \item [{[13.]}] Varadarajan, V.S., \textit{``An introduction to harmonic analysis on semisimple Lie groups,''}  Cambridge University Press, Cambridge, $1989.$
                        \end{description}

\end{document}